\documentclass{amsart}
\usepackage{amssymb,latexsym,color}
\usepackage{graphicx,pgf,pst-all}
\theoremstyle{plain}
\newtheorem{theorem}{Theorem}

\newtheorem{lemma}{Lemma}

\theoremstyle{definition}

\theoremstyle{remark}

\newtheorem{remark}{\bf Remark}
\numberwithin{equation}{section}
\newcommand{\pmx}{p_{\text{max}}}
\newcommand{\pmn}{p_{\text{min}}}
\newcommand{\bp}{\bar{p}}
\newcommand{\ubknp}{(u-\bar{k}_n)_+}
\newcommand{\uknp}{(u-{k}_n)_+}
\newcommand{\ine}{\int_{E}}

\newcommand{\noi}{\noindent}

\newcommand{\dsty}{\displaystyle}
\newcommand{\txty}{\textstyle}

\newcommand{\al}{\alpha}

\newcommand{\gm}{\gamma}

\newcommand{\kp}{\kappa}
\newcommand{\varep}{\varepsilon}

\newcommand{\vp}{\varphi}
\newcommand{\sig}{\sigma}

\newcommand{\om}{\omega}

\newcommand{\z}{\zeta}

\newcommand{\df}[1]{\buildrel\mbox{\small def}\over{#1}}

\newcommand{\nn}{\mathbb{N}}

\newcommand{\rr}{\mathbb{R}}
\newcommand{\rn}{\rr^N}

\newcommand{\bl}[1]{\mathbf{#1}}

\newcommand{\dvg}{\operatorname{div}}

\newcommand{\essup}{\operatornamewithlimits{ess\,sup}}
\newcommand{\essinf}{\operatornamewithlimits{ess\,inf}}
\newcommand{\essosc}{\operatornamewithlimits{ess\,osc}}

\newcommand{\loc}{\operatorname{loc}}

\newcommand{\dist}{\operatorname{dist}}

\newcommand{\ttop}[2]{\genfrac{}{}{0pt}{}{#1}{#2}}

\newcommand{\bsu}{\mathop{\txty{\sum}}\limits}
\newcommand{\pro}{\mathop{\txty{\prod}}\limits}

\newcommand{\pl}{\partial}

\def\Xint#1{\mathchoice
    {\XXint\displaystyle\textstyle{#1}}%
    {\XXint\textstyle\scriptstyle{#1}}%
    {\XXint\scriptstyle\scriptscriptstyle{#1}}%
    {\XXint\scriptscriptstyle\scriptscriptstyle{#1}}%
    \!\int}
\def\XXint#1#2#3{\setbox0=\hbox{$#1{#2#3}{\int}$}
    \vcenter{\hbox{$#2#3$}}\kern-0.5\wd0}
\def\bint{\Xint-}
\def\dashint{\Xint{\raise4pt\hbox to7pt{\hrulefill}}}
\def\dashiint{\bint\kern-0.15cm\bint}

\newcommand{\cio}[1]{C_o^\infty\!\left({#1}\right)}

\begin{document}
\title[Anisotropic $p$-Laplacian Equation]
      {Remarks on Local Boundedness and Local H\"older 
Continuity of Local Weak Solutions to Anisotropic 
$p$-Laplacian Type Equations}
      
\author{Emmanuele DiBenedetto}
\address{Department of Mathematics\\
               Vanderbilt University\\  
               1326 Stevenson Center\\
               Nashville TN 37240, USA}
\email{em.diben@vanderbilt.edu}
\urladdr{http://www.math.vanderbilt.edu/~dibe/}
\author{Ugo Gianazza}
\address{Dipartimento di Matematica ``F. Casorati"\\
               Universit\`a di Pavia\\ 
               via Ferrata 1\\ 
               I-27100 Pavia}
\email{gianazza@imati.cnr.it}
\urladdr{http://arturo.imati.cnr.it/~gianazza}
\author{Vincenzo Vespri}
\address{Dipartimento di Matematica e Informatica ``U. Dini"\\ 
               Universit\`a di Firenze\\ 
               viale Morgagni 67/A, 
               I-50134 Firenze}
\email{vespri@math.unifi.it}
\urladdr{http://web.math.unifi.it/users/vespri/}
      
\thanks{E. DiBenedetto's research supported  under NSF grant DMS-1265548}  
\keywords{Anisotropic $p$-Laplacian, Elliptic, H\"older continuity} 
\subjclass[2010]{Primary: 35J70, 35J92, 35B65; 
Secondary 35B45}
\date{}
\begin{abstract}
Locally bounded, local weak solutions to a special class 
of quasilinear, ani\-so\-tropic, $p$-Laplacian type elliptic 
equations, are shown to be locally H\"older continuous. 
Homogeneous local upper bounds are established for local 
weak solutions to a general class of quasilinear 
anisotropic  equations.
\end{abstract}
\maketitle
\section{Introduction}
Consider quasi-linear, elliptic differential equations of the form
\begin{equation}\label{equation}
\dvg\bl{A}(x,u,Du) = 0\quad
\text{ weakly in sone open set }\> E\subset \rn
\end{equation}
where the function $\bl{A}=(A_1,\dots,A_N):E\times\rr^{N+1}\to\rn$ 
is only assumed to be measurable
and subject to the structure conditions 
\begin{equation}\label{equation1}
\begin{aligned}
A_{i}(x,u,Du)\cdot u_{x_i}&\ge C_{o,i} |u_{x_i}|^{p_i},\\
|A_{i}(x,u,Du)|&\le C_{1,i}|u_{x_i}|^{p_i-1},\\
\end{aligned}%
\end{equation}
where $p_i>1$ and $C_{o,i}$ and $C_{1,i}$ are given positive 
constants.  Such elliptic equations are termed \emph{anisotropic}, 
their prototype being
\begin{equation}\label{equation-proto}
\bsu_{i=1}^N\big(\left|u_{x_i}\right|^{p_i-2}u_{x_i}\big)_{x_i}= 0 \quad
\text{ in }\ E.
\end{equation}
For a multi-index ${\bf p}=\{p_1,\dots,p_N\}$, $p_i\ge1$, let
\[ 
W^{1,{\bf p}}({E})=\{u\in L^{1}({E}):\ u_{x_i}\in L^{p_i}({E}),\ i=1,\dots,N\},
\]
and
\[
W_o^{1,{\bf p}}({E})=W^{1,{\bf p}}({E})\cap W_o^{1,1}({E}).
\]
A function $\dsty u\in W^{1,{\bf p}}_{\loc}({E})$
is a local, weak solution to \eqref{equation} if for every compact 
set $K\subset E$
\begin{equation}\label{Eq:1:4}
\int_K\bl{A}(x,u,Du)\cdot D\vp\,dx=0\quad\text{ for all }\> \vp\in\cio{K}.
\end{equation}
The parameters $\{N, p_i, C_{o,i}, C_{1,i}\}$ are the data, and 
we say that a generic constant $\gm=\gm(N,p_i,C_{o,i},C_{1,i})$ 
depends upon the data, if it can be quantitatively determined 
a priori only in terms of the indicated parameters.

Define, 
\begin{equation}\label{Eq:1:5}
\frac{1}{\bp}=\frac{1}{N}\bsu_{i=1}^N\frac{1}{p_{i}},\qquad\qquad
\begin{array}{c}
\pmn=\min\{p_1,\dots,p_N\},\\
\pmx=\max\{p_1,\dots,p_N\}.
\end{array}
\end{equation}
For a compact set $K\subset E$ introduce the {\it intrinsic}, elliptic
${\bf p}$-distance from $K$ to $\partial E$ by 
\begin{equation*}
{\bf p}-\dist(K;\partial E)\,\df{=}\,
\inf_{\ttop{x\in K}{y\in\partial E}}
{\left(\bsu_{j=1}^N\|u\|_{\infty,E}^{\frac{\pmx-p_j}{\pmx}} |x_j-y_j|^{\frac{p_j}{\pmx}}\right)}. 
\end{equation*} 
\begin{theorem}\label{Thm:1:1}
Let $u$ be a bounded, local, weak solution to 
\eqref{equation}--\eqref{equation1}, and assume $\bar p<N$. 
There exists a positive quantity $q>1$, depending only on the 
data, such that if 
\begin{equation}\label{Eq:smallness}
p_{\rm max}-p_{\rm min}\le\frac1q,
\end{equation}
then $u$ is locally H\"older continuous in $E$, i.e. there exist 
constants $\gm>1$ and $\al\in(0,1)$ depending only on the data, 
such that for every compact set $K\subset E$,
\begin{equation*}
|u(x_1)-u(x_2)|\le\gm\|u\|_{\infty,E}
{\Bigg(\frac{\bsu_{i=1}^N\|u\|_{\infty,E}^{\frac{p_{\rm max}-p_i}{p_{\rm max}}}
|x_{1,i}-x_{2,i}|^{\frac{p_i}{p_{\rm max}}}}{{\bf p}-\dist(K;\partial E)}\Bigg)^\al}
\end{equation*}
for every pair of points $x_1,\,x_2\,\in K$.
\end{theorem}
\begin{remark}{\normalfont
For a general distribution of the $p_j$, 
unbounded weak solutions might exist (\cite{G,M0}). 
In \cite{FS,K,CMM} it is shown that local weak solutions are locally
bounded provided
\begin{equation}\label{Eq:boundedness-condition}
\bar p<N,\qquad\pmx\le\frac{N\bar p}{N-\bar p}.
\end{equation} 
In Section~\ref{S:bound} we revisit and improve these boundedness 
estimates.
}
\end{remark}
\begin{remark}\label{Rmk:holder-intro}
{\normalfont 
The H\"older continuity of solutions holds also for $\bar p\ge N$; indeed, when $\bar p=N$, a straightforward modification of our arguments suffices, strictly analogous to the one used in the isotropic case when $p=N$; when $\bar p>N$, the embedding results of \cite[Theorem~2]{HaSch} (se also \cite{Ra}), ensure that $u$ is H\"older continuous.}
\end{remark}
\begin{remark}\label{Rmk:1:1} {\normalfont
The constants $\gm$ and $\al$ deteriorate as 
either $p_i\to\infty$ or $p_i\to1$, in the sense that 
$\gm(p)\to\infty$ and $\al\to0$ as either $p_i\to\infty$, or $p_i\to1$.
}
\end{remark}
\vskip.2truecm

\noi{\it Acknowledgements - } We thank Prof. P. Marcellini, 
E. Mascolo and G. Cupini for enlightening conversations on 
embeddings for anisotropic Sobolev spaces.
\section{Novelty and Significance}
If the coefficients in \eqref{equation}--\eqref{equation1} are 
differentiable, and satisfy some further, suitable structure 
conditions, Lipschitz estimates have been derived by 
Marcellini \cite{M,M1}. 
If the coefficients are merely bounded and measurable,
H\"older continuity has been established in \cite{LS} in 
the special case of $p_1=2<p_2=p_3=\cdots=p_N$, i.e., 
the $p_j$ are all the same except the smallest one. 
The main idea is to regard the equation as ``parabolic'' with 
respect to the variable $x_1$, corresponding to $p_1$, and 
to apply the techniques of \cite{DB93,DBGV}. An extension 
to the case $1<p_1<p_2=p_3=\cdots=p_N$, by the 
same techniques, is in \cite{DMV}. 

Theorem~\ref{Thm:1:1} is a further step in understanding the 
regularity of solutions of anisotropic elliptic equations, 
with full quasi-linear structure. Our approach 
is ``elliptic'' in nature, it is modelled after \cite{degiorgi}, 
no variable is regarded as ``parabolic'', and no restriction 
is placed on the distribution of the $p_j$ other that $\pmx-\pmn\ll1$. 
In particular, the $p_j$ could all be different. 

While partial, Theorem~\ref{Thm:1:1} disproves the claim 
in \cite{G} by which, H\"older continuity of weak solutions 
to \eqref{equation}--\eqref{equation1}, holds if and only if 
$\pmn=\pmx$, i.e., if no anisotropy is present. 

Finally, Theorem~\ref{Thm:1:1} can be seen as a stability result of 
the H\"older continuity of solutions, when $p_i\to\pmn$, and 
correspondingly, the anisotropic $p$-laplacian tends to 
the $\pmn$-laplacian.
\section{Preliminaries and Intrinsic Geometry}
\begin{lemma}[Sobolev-Troisi Inequality, \cite{troisi}]\label{troisi}
Let ${E}\subset\rn$ be a bounded, open set and consider $u\in
W_o^{1,{\bf p}}({E}),$ $p_{i}>1$  for all
$i=1,\dots,N$. Assume $\bar{p}<N$ and let
\begin{equation}\label{Eq:S:intr:16}
{p}_{*}=\frac {N\bar{p}}{N-\bar{p}}.
\end{equation}
Then there exists a constant $c$ depending only on $N,{\bf p}$,
such that
\begin{equation*}
\left\Vert u \right\Vert_{L^{{p}_{*}}\left({E}\right)
}^{N}\leq c\,\pro_{i=1}^N\left\Vert u_{x_{i}}
\right\Vert_{L^{p_{i}}\left({E}\right)}.
\end{equation*}
\end{lemma}
For $\rho>0$ consider the cube $K_\rho=(-\rho,\rho)^N$, with 
center at the origin of $\rn$ and edge $2\rho$,  and set
\begin{equation}\label{Eq:S:intr:1}
\mu^+=\essup_{K_{2\rho}} u;\qquad
\mu^-=\essinf_{K_{2\rho}} u;\qquad
\om=\mu^+-\mu^-=\essosc_{K_{2\rho}} u.
\end{equation}
These numbers being determined, construct the cylinder
\begin{equation}\label{Eq:S:intr:2}
Q_{\rho}=\pro_{j=1}^N(-\rho_j,\rho_j),
\end{equation}
with $0<\rho_j\le\rho$ to be determined.
This implies that $Q_{2\rho}\subset K_{2\rho}$ and hence
$\essosc_{Q_{\rho}} u\le \om$. 
\subsection{Basic Equation and Energy Inequalities}\label{S:intr:1}
For $\sig\in(0,1)$ let $\z_j$ be a non-negative, piecewise 
smooth cutoff function in the interval $(-\rho_j,\rho_j)$ 
which equals 1 on $(-\sig\rho_j,\sig\rho_j)$, vanishes at 
$\pm\rho_j$, and such that $|\z'_{j}|\le[(1-\sig)\rho_j]^{-1}$. 

Set $\z={\pro}_{j=1}^N\z_j^{p_j}$, and in the weak formulation of
\eqref{equation}--\eqref{equation1}
take the testing function $\pm(u-k)_\pm\z$. This gives, after 
standard calculations, 
\begin{equation}\label{Eq:S:intr:7}
\bsu_{j=1}^N\int_{Q_\rho}
\Big|\frac{\pl}{\pl x_j}\big[(u-k)_{\pm}\z^{\frac1{p_j}}]\Big|^{p_j} dx \le 
\gm\bsu_{j=1}^N\frac1{(1-\sig)^{p_j}\rho_j^{p_j}}
\int_{Q_\rho}(u-k)_{\pm}^{p_j}dx. 
\end{equation}
The constant $\gamma$ depends only upon the 
data, and is independent of 
$\rho$.
\section{DeGiorgi Type Lemmas}\label{S:intr:2}
Taking $k=\mu^+-\frac{\om}{2^s}$, for $s\ge1$, 
and $(u-k)_+$ in (\ref{Eq:S:intr:7}) yields
\begin{equation}\label{Eq:S:intr:8}
\begin{aligned}
&\bsu_{j=1}^N\int_{Q_\rho} 
\Big|\Big[\Big(u-\Big(\mu^+-\frac{\om}{2^s}\Big)\Big)_+
\z^{\frac1{p_j}}\Big]_{x_j}\Big|^{p_j}dx\\
&\le \frac{\gm}{(1-\sig)^{\pmx}}
\bsu_{j=1}^N\frac1{\rho_j^{p_j}}\Big(\frac{\om}{2^s}\Big)^{p_j} 
\big|Q_\rho\cap[u>\mu^+-\frac{\om}{2^s}]\big|. 
\end{aligned}
\end{equation}
Likewise, taking $k=\mu^-+\frac{\om}{2^s}$, for $s\ge1$, 
and $-(u-k)_-$ in (\ref{Eq:S:intr:7}) yields
\begin{equation}\label{Eq:S:intr:9}
\begin{aligned}
&\bsu_{j=1}^N\int_{Q_\rho} 
\Big|\Big[\big(u-\Big(\mu^-+\frac{\om}{2^s}\Big)\Big)_-
\z^{\frac1{p_j}}\Big]_{x_j}\Big|^{p_j}dx\\
&\le \frac{\gm}{(1-\sig)^{\pmx}}
\bsu_{j=1}^N\frac1{\rho_j^{p_j}}\Big(\frac{\om}{2^s}\Big)^{p_j} 
\big|Q_\rho\cap[u<\mu^-+\frac{\om}{2^s}]\big|. 
\end{aligned}
\end{equation}
Choose 
\begin{equation}\label{Eq:intr:radius}
\rho_j=\left(\frac\om{2^q}\right)\rho^{\frac\al{p_j}}\quad\text{
 for some $q>0$ and $\al\ge\pmx$ to be chosen} 
\end{equation}
and let $Q_\rho$ the cylinder in \eqref{Eq:S:intr:2} for such 
a choice of $\rho_j$.  Without loss of generality, may assume 
$\om\le1$, so that $0<\rho_j\le\rho$ as required.  
\begin{lemma}\label{Lm:S:intr:1}
There exists a number $\nu\in(0,1)$ depending only upon the data, 
such that if 
\begin{equation}\label{Eq:S:intr:10piu}
\Big|\Big[u>\mu^+-\frac{\om}{2^q}\Big]\cap Q_\rho
\Big |<\nu \big|Q_\rho\big|,
\end{equation}
for some $q\in\nn$, then 
\begin{equation}\label{Eq:S:intr:11piu}
u\le\mu^+-\frac{\om}{2^{q+1}}\quad\text{ a.e. in }\quad Q_{\frac12\rho}.
\end{equation}
\end{lemma}
Likewise
\begin{lemma}\label{Lm:S:intr:2}
There exists a number $\nu\in(0,1)$ depending only upon the data, 
such that if 
\begin{equation}\label{Eq:S:intr:10me}
\Big|\Big[u<\mu^-+\frac{\om}{2^q}\Big]\cap Q_\rho
\Big |<\nu \big|Q_\rho\big|,
\end{equation}
for some $q\in\nn$, then 
\begin{equation}\label{Eq:S:intr:11me}
u\ge\mu^-+\frac{\om}{2^{q+1}}\quad\text{ a.e. in }
\quad Q_{\frac12\rho}.
\end{equation}
\end{lemma}
We prove only Lemma~\ref{Lm:S:intr:1}, the proof of 
Lemma~\ref{Lm:S:intr:2} being analogous.
\vskip0.3cm
\begin{proof}  

For each $j\in\{1,\dots,N\}$ consider the sequence of radii
\begin{equation}\label{Eq:S:intr:11}
\rho_{j,n}=\frac12\rho_j\Big(1+\frac1{2^n}\Big),\qquad 
\text{ for }\quad n=0,1,\dots.
\end{equation}
This is a decreasing sequence with $\rho_{j,o}=\rho_j$ 
and $\rho_{j,\infty}=\frac12\rho_j$. The corresponding cylinders
\begin{equation}\label{Eq:S:intr:12}
Q_n\df{=}Q_{\rho_n}={\pro}_{j=1}^N\big(-\rho_{j,n},\rho_{j,n}\big)
\qquad \text{ for }\quad n=0,1,\dots
\end{equation}
are nested, i.e., $Q_{n+1}\subset Q_n$, with $Q_o=Q_\rho$ 
and $Q_{\infty}= Q_{\frac12\rho}$, since $\al\ge1$.  
For each $j\in\{1,\dots,N\}$ let $\z_{j,n}$ be a standard 
non-negative cutoff function in $(-\rho_{j,n},\rho_{j,n})$ which 
equals 1 on $(-\rho_{j,n+1},\rho_{j,n+1})$, vanishes 
at $\pm\rho_{j,n}$ and such that 
$|\z^\prime_{j,n}|\le2^{n+2}\rho_{j,n}^{-1}$. Then set 
$\z_n={\pro}_{j=1}^N\z_{j,n}^{p_j}$ to be a cutoff function in 
$Q_n$ that equals 1 on $Q_{n+1}$.  Consider also the increasing 
sequence of levels
\begin{equation}\label{Eq:S:intr:13}
k_n=\mu^+-\frac{\om}{2^{q+1}}\Big(1+\frac1{2^n}\Big),\qquad 
\text{ for }\quad n=0,1,\dots
\end{equation}
and in the weak formulation of \eqref{equation}--\eqref{equation1}, 
take the test function 
$(u-k_n)_+\z_n$. This leads to analogues of 
\eqref{Eq:S:intr:8} over the cylinders $Q_n$, 
with $1-\sig>2^{-(n+2)}$, and $q\le s<q+1$. Rewriting \eqref{Eq:S:intr:8} 
with these specifications gives
\begin{equation}\label{Eq:S:intr:14}
\bsu_{j=1}^N\int_{Q_n}\big|
\big[(u-k_n)_+\z_n^{\frac1{p_j}}\big]_{x_j}\big|^{p_j}dx
\le\gm \frac{2^{n\pmx}}{\rho^{\al}}\big|Q_n\cap[u>k_n]\big|.
\end{equation}
Since $(u-k_n)_+\z_n$ vanishes on $\pl Q_n$, by the
anisotropic embedding of Lemma~\ref{troisi} 
\begin{equation}\label{Eq:S:intr:15}
\Big(\int_{Q_n}\big[(u-k_n)_+\z_n\big]^{p_*} dx
\Big)^{\frac1{p_*}}\le c\,\pro_{j=1}^N\Big(\int_{Q_n}\big|\big[
(u-k_n)_+\z_n\big]_{x_j}\big|^{p_j}dx\Big)^{\frac1{N p_j}}.
\end{equation}
where $p_*$ is as in \eqref{Eq:S:intr:16}, $p_j>1$  for 
$j=1,\dots,N$, $\bar p<N$. Since $0\le\z_n\le1$ and $p_j>1$ estimate
\begin{equation*}
\int_{Q_n}\big|\big[(u-k_n)_+\z_n\big]_{x_j}\big|^{p_j}dx\le 
\gm \int_{Q_n}\big|\big[(u-k_n)_+\z_n^{\frac1{p_j}}\big]_{x_j}\big|^{p_j}dx.
\end{equation*}
Therefore, combining this with (\ref{Eq:S:intr:15}) and 
(\ref{Eq:S:intr:14}) gives
\begin{equation*}
\begin{aligned}
&(k_{n+1}-k_n)\big|Q_{n+1}\cap [u>k_{n+1}]\big|\le 
\int_{Q_{n+1}\cap [u>k_{n+1}]}(u-k_n)_+dx\\
&\le \int_{Q_n}(u-k_n)_+\z_n dx\le\Big(\int_{Q_n}\big[(u-k_n)_+
\z_n\big]^{p_*} dx\Big)^{\frac1{p_*}} 
\big|Q_n\cap[u>k_n]\big|^{1-\frac1{p_*}}\\
&\le\gm\pro_{j=1}^N\Big(\int_{Q_n}\big|\big[
(u-k_n)_+\z_n^{\frac1{p_j}}\big]_{x_j}\big|^{p_j}dx\Big)^{\frac1{N p_j}}.
\big|Q_n\cap[u>k_n]\big|^{1-\frac1{p_*}}\\
&\le\gm \Big(2^{\frac{\pmx}{\bp}}\Big)^n 
\frac1{\dsty \rho^{\frac{\al}{\bp}}}
\big|Q_n\cap[u>k_n]\big|^{1+\frac1{\bp}-\frac1{p_*}}\\
&=\gm b^n \frac{\dsty\Big(\frac{\om}{2^q}\Big)}{\big|Q_\rho\big|^{\frac1N}}
\big|Q_n\cap[u>k_n]\big|^{1+\frac1N}
\end{aligned}
\end{equation*}
where we have set $b=2^{\pmx/\bp}$.
By the definition of $k_n$ in (\ref{Eq:S:intr:13}), 
the first term in round brackets on the left hand side is 
\begin{equation*}
k_n-k_{n+1}=\Big(\frac{\om}{2^q}\Big)\frac1{2^{2+n}}.
\end{equation*}
Combining these remarks and inequalities, and setting
\begin{equation*}
Y_n=\frac{\big|Q_n\cap[u>k_n]\big|}{|Q_\rho|}
\end{equation*}
yields the recursive inequalities
\begin{equation*}
Y_{n+1}\le C(2b)^nY_n^{1+\frac1N}
\end{equation*}
for constants $C$ and $b$ depending only upon the data. 
It follows from these and Lemma~5.1 of \cite[Chapter~2]{DBGV}, 
that there exists a number $\nu\in(0,1)$ depending only 
on $\{C,b,N\}$, and hence only upon the data,
such that  $\{Y_n\}\to0$ as $n\to\infty$ provided
\begin{equation*}
Y_o=\frac{\big|Q_\rho\cap\big[u>\mu^+-\frac{\om}{2^q}
\big]\big|}{|Q_\rho|}\le\nu.
\end{equation*}
\end{proof}
By these lemmata the H\"older continuity of $u$ will follow
by standard arguments, if one can determine $q$, and hence the 
intrinsic cylinders $Q_\rho$, for which either 
\eqref{Eq:S:intr:10piu} or \eqref{Eq:S:intr:10me} holds.
\section{Proof of Theorem~\ref{Thm:1:1}}\label{S:small} 
Assume that 
\begin{equation}\label{Eq:small:1}
\big|\big[u<\mu^-+{\txty\frac12}\om\big]\cap Q_\rho\big|\ge{\txty\frac12}
\big|Q_\rho\big|.
\end{equation}
For each $s\in\nn$ with $s\le q$, introduce the two complementary sets
\begin{equation}\label{Eq:small:2}
A_s=\Big[u>\mu^+-\frac{\om}{2^s}\Big]\cap Q_\rho;\qquad
Q_\rho-A_s=\Big[u\le\mu^+-\frac{\om}{2^s}\Big]\cap Q_\rho
\end{equation}
and consider the doubly truncated function 
\begin{equation}\label{Eq:small:3}
 v_s=\left\{
\begin{array}{ll}
0\quad&{\dsty \text{for }\quad\kern1.8cm u<\mu^+-\frac{\om}{2^s}},\\
{}\\
{\dsty u-\Big(\mu^+-\frac{\om}{2^s}\Big)}\quad&{\dsty \text{for }\quad
\mu^+-\frac{\om}{2^s}\le u <\mu^+-\frac{\om}{2^{s+1}}},\\
{}\\
{\dsty \frac{\om}{2^{s+1}}}\quad&{\dsty \text{for }\quad
 \mu^+-\frac{\om}{2^{s+1}}\le u.}
\end{array}\right.
\end{equation}
By construction $v_s$ vanishes on $Q_\rho-A_s$. Pick any two points 
\begin{equation*}
x=(x_1,\dots,x_N)\in A_s\quad\text{ and }\quad 
y=(y_1,\dots,y_N)\in Q_\rho-A_s
\end{equation*}
and construct a polygonal joining $x$ and $y$ and sides parallel 
to the coordinate axes, say for example $P_N=x$ and
\begin{align*}
P_{N-1}&=(x_1,\dots,x_{N-1},y_N);\quad 
P_{N-2}=(x_1,x_2,\dots,y_{N-1},y_N);\quad \cdots\quad\\ 
P_1&=(x_1,y_2,\dots,y_N);\quad P_o=(y_1,\dots,y_N).
\end{align*}
By elementary calculus
\begin{equation*}
\begin{aligned}
v_s(x)&= \left[v_s(P_N)-v_s(P_{N-1})\right]+\cdots+
\left[v_s(P_1)-v_s(P_o)\right]\\
&=\int_{y_N}^{x_N}\frac{\pl}{\pl x_N}v_s(x_1,\dots,x_{N-1},t)dt+ 
\int_{y_{N-1}}^{x_{N-1}} \frac{\pl}{\pl x_{N-1}}
v_s(x_1,\dots,x_{N-2},t,y_N)dt\\
&\quad+\dots+\int_{y_1}^{x_1}\frac{\pl}{\pl x_1}
v_s(t,y_2,\dots,y_N)dt\\
&\le \bsu_{j=1}^N\int_{-\rho_j}^{\rho_j} 
\big|v_{s,x_j}\big|(x_1,\dots,
\underbrace{t}_{j-\text{th}\ variable},\dots,y_N)dt
\end{aligned}
\end{equation*}
where the quantities $\rho_j$ are defined in (\ref{Eq:intr:radius}). 
Integrate in $dx$ over 
$A_s$ and in $dy$ over $Q_\rho-A_s$, and 
take into account (\ref{Eq:small:1}) to get
\begin{equation*}
\frac12\big|Q_\rho\big|\int_{Q_\rho} v_s dx
\le2\big|Q_\rho\big|
\bsu_{j=1}^N\rho_j \int_{Q_\rho}
\big|v_{s,x_j}\big|dx.
\end{equation*}
From this, by the definitions (\ref{Eq:small:2}) 
and (\ref{Eq:small:3}) of $A_s$ and $v_s$,
\begin{equation}\label{Eq:small:4}
\begin{aligned}
\frac{\om}{2^{s+1}}\big|A_{s+1}\big|&\le 
4\bsu_{j=1}^N\rho_j\int_{A_s-A_{s+1}}
\big|u_{x_j}\big|dx\\
&\le4\bsu_{j=1}^N\rho_j\left(\int_{A_s-A_{s+1}}
\big|u_{x_j}\big|^{\pmn}dx\right)^{\frac1{\pmn}} 
\big|A_s-A_{s+1}\big|^{1-\frac1{\pmn}}\\
&\le4\bsu_{j=1}^N\rho_j\left(\int_{A_s-A_{s+1}}
\big|u_{x_j}\big|^{p_j}dx\right)^{\frac1{p_j}}
\big|Q_\rho\big|^{\frac1{\pmn}-\frac1{p_j}} 
\big|A_s-A_{s+1}\big|^{1-\frac1{\pmn}}.
\end{aligned}
\end{equation}
For each $j$ fixed, the integrals involving $u_{x_j}$ are 
estimated by means of \eqref{Eq:S:intr:8} applied over the 
pair of cubes $Q_\rho$ and $Q_{2\rho}$, as follows:
\begin{equation*}
\begin{aligned}
\left(\int_{A_s-A_{s+1}}\big|u_{x_j}\big|^{p_j}dx
\right)^{\frac1{p_j}}&\le\left(\int_{Q_{\rho}}\left|\frac{\pl}{\pl x_j} 
\left(u-\left(\mu^+-\frac{\om}{2^s}\right)\right)_+\right|^{p_j}dx
\right)^{\frac1{p_j}}\\
&\le \gm \left(\frac1{\rho^{\al}}\bsu_{\ell=1}^N\left(\frac{\om}{2^s}
\right)^{p_\ell}\left(\frac{\om}{2^q}\right)^{-p_\ell}
\big|Q_{\rho}\big| \right)^{\frac1{p_j}}\\
&\le \gm \left(\frac1{\rho^{\al}}\bsu_{\ell=1}^N
\left(\frac{2^q}{2^s}\right)^{p_\ell}
\big|Q_{\rho}\big| \right)^{\frac1{p_j}}\\
&=\gm\left(\frac1{\rho_j^{p_j}}
\left(\frac{\om}{2^q}\right)^{p_j}\bsu_{\ell=1}^N
\left(\frac{2^q}{2^s}\right)^{p_\ell}
\big|Q_{\rho}\big|\right)^{\frac1{p_j}}.
\end{aligned}
\end{equation*}
If $p_\ell\le p_j$, since $s\le q$ estimate
\begin{equation}\label{Eq:small:5}
\left(\frac{2^q}{2^s}\right)^{p_\ell}
\le\left(\frac{2^q}{2^s}\right)^{p_j}= 
\left(\frac{\om}{2^s}\right)^{p_j} 
\left(\frac{\om}{2^q}\right)^{-p_j},\qquad(\text{ case of }\> p_\ell\le p_j). 
\end{equation}
If $p_\ell>p_j$ since $s\le q$ compute and estimate 
\begin{equation}\label{Eq:small:6}
\begin{aligned}
\left(\frac{2^q}{2^s}\right)^{p_\ell}
&=\left(\frac{2^q}{2^s}\right)^{p_j} 
\left(\frac{2^q}{2^s}\right)^{p_\ell-p_j}\\
&\le\left(\frac{\om}{2^s}\right)^{p_j} 
\left(\frac{\om}{2^q}\right)^{-p_j}2^{q(\pmx-\pmn)}
\end{aligned}
\qquad(\text{case of }\> p_\ell>p_j). 
\end{equation}
Assume momentarily that the number $q$ has been chosen. 
Then stipulate that $q(\pmx-\pmn)\le1$. For such a choice 
we have in all cases 
\begin{equation*}
\left(\int_{A_s-A_{s+1}}\big|u_{x_j}\big|^{p_j}dx
\right)^{\frac1{p_j}}\le\gm \frac1{\rho_j}
\left(\frac{\om}{2^s}\right) 
\big|Q_\rho\big|^{\frac1{p_j}}.  
\end{equation*}
Combining these estimates in (\ref{Eq:small:4}) yields
\begin{equation*}
\big|A_{s+1}\big|\le\gm \big|Q_{\rho}\big|^{\frac1{\pmn}} 
\big(\big|A_s\big|-\big|A_{s+1}\big|\Big)^{1-\frac1{\pmn}}.
\end{equation*}
Take the $\big(\frac{\pmn}{{\pmn}-1}\big)$ power and add for 
$s=1,\dots (q-1)$ to get
\begin{equation*}
(q-1)\big|A_q\big|^{\frac{\pmn}{\pmn-1}}
\le\gm^{\frac{\pmn}{\pmn-1}}\big|Q_{\rho}\big|^{\frac1{\pmn-1}} 
\big|A_o\big|.
\end{equation*}
From this 
\begin{equation*}
\big|A_q\big|\le\frac{\gm}{(q-1)^{\frac{\pmn-1}{\pmn}}}
\big|Q_{\rho}\big|.
\end{equation*}
In the DeGiorgi-type Lemma, the number $\nu$ is 
independent of $q$. Now choose $q$ so that 
\begin{equation}\label{Eq:small:7}
\big|A_q\big|\le\nu \big|Q_{\rho}\big|,\quad\text{ for }\quad
\nu=\frac{\gm}{(q-1)^{\frac{\pmn-1}{\pmn}}}.
\end{equation}
Notice that $q$ is determined in terms of $\pmn$ and not 
in terms of the difference $(\pmx-\pmn)$. Thus, one determines 
first $q$ from (\ref{Eq:small:7}) in terms only of the data. 
Then (\ref{Eq:smallness}), for such a choice of $q$, serves 
as a condition of H\"older continuity for $u$. 
\section{Boundedness}\label{S:bound}
Continue to denote by $u\in W^{1,\bl{p}}_{\loc}(E)$ a 
local weak solution to \eqref{equation}--\eqref{equation1}, 
in the sense of (\ref{Eq:1:4}), {with $\bar p<N$}. The estimations below 
use that $u\in L^{p_*}_{\loc}(E)$. When $\pmx<p_*$ this is 
insured by the embeddings in \cite[Theorem~1]{KK} or 
\cite{nikolski}. If $\pmx=p_*$ in what follows the membership 
$u\in L^{p_*}_{\loc}(E)$, is assumed.

When $\bar p=N$, then \cite[Theorem~1]{HaSch} ensures that $u\in L^q_{\loc}(E)$ for any arbitrary $1\le q<\infty$, and the arguments below can be repeated verbatim, to obtain a quantitative estimate of the local boundedness of $u$. Finally, when $\bar p>N$, as we mentioned in Remark~\ref{Rmk:holder-intro}, proper Morrey-type embeddings directly ensure the boundedness of $u$.
\subsection{Some General Recursive Inequalities}
Let $\rho_j$ as in (\ref{Eq:S:intr:2}) to be defined, and for each 
$j\in\{1,\dots,N\}$ consider the sequence of radii, 
\begin{equation}\label{Eq:S:bdd:5}
\rho_{j,n}=\frac12\rho_j\Big(1+\frac1{2^n}\Big),\qquad 
\text{ for }\quad n=0,1,\dots.
\end{equation}
This is a decreasing sequence, with $\rho_{j,n}=\rho_j$ 
and $\rho_{j,\infty}=\frac12\rho_j$. The corresponding cylinders
\begin{equation}\label{Eq:S:bdd:6}
Q_n\df{=}Q_{\rho_n}={\pro}_{j=1}^N\big(-\rho_{j,n},\rho_{j,n}\big)
\qquad \text{ for }\quad n=0,1,\dots
\end{equation}
are nested, i.e., $Q_{n+1}\subset Q_n$, with $Q_o=Q_\rho$ 
and $Q_{\infty}= Q_{\frac12\rho}$, since $\al\ge p_j$.  
For each $j\in\{1,\dots,N\}$ let $\z_{j,n}$ be a standard 
non-negative cutoff function in $(-\rho_{j,n},\rho_{j,n})$ which 
equals 1 on $(-\rho_{j,n+1},\rho_{j,n+1})$, vanishes 
at $\pm\rho_{j,n}$ and such that 
$|\z^\prime_{j,n}|\le2^{n+2}\rho_{j,n}^{-1}$. Then set 
$\z_n={\pro}_{j=1}^N\z_{j,n}^{p_j}$ to be a cutoff function in 
$Q_n$ that equals 1 on $Q_{n+1}$.

Consider also the increasing 
sequence of levels
\begin{equation}\label{Eq:S:bdd:7}
k_n=\Big(1-\frac1{2^n}\Big) k,\quad\text{ and }\quad
\bar{k}_n=\frac{k_{n+1}+k_n}2
\end{equation}
for $n=0,1,\dots$, with $k>0$ to be chosen. By the definition 
$k_o=0$ and $k_\infty=k$.  Write 
(\ref{Eq:S:intr:7}) for $\ubknp\z_n$, over the 
cylinders $Q_n$, with $1-\sig>2^{-(n+2)}$. 
Since $\ubknp\z_n$ vanishes on $\pl Q_n$, by the
anisotropic embedding of Lemma~\ref{troisi}
\begin{equation}\label{Eq:S:bdd:8}
\Big(\int_{Q_n}\big[\ubknp\z_n\big]^{p_*} dx
\Big)^{\frac1{p_*}}\le\gamma\pro_{j=1}^N\Big(\int_{Q_n}\big|\big[
\ubknp\z_n\big]_{x_j}\big|^{p_j}dx\Big)^{\frac1{N p_j}}.
\end{equation}
where $p_*$ has been defined in \eqref{Eq:S:intr:16}.

Since $0\le\z_n\le1$ and $p_j\ge1$ estimate
\begin{equation*}
\int_{Q_n}\big|\big[\ubknp\z_n\big]_{x_j}\big|^{p_j}dx\le 
\gm \int_{Q_n}\big|\big[\ubknp\z_n^{\frac1{p_j}}\big]_{x_j}\big|^{p_j}dx.
\end{equation*}
Therefore, combining this with (\ref{Eq:S:intr:7}) and 
(\ref{Eq:S:bdd:8}) gives
\begin{equation*}
\begin{aligned}
\left(\int_{Q_{n+1}}\ubknp^{p_*} dx\right)^{\frac1{p_*}}&\le
\left(\int_{Q_n}\big[\ubknp\z_n\big]^{p_*} dx\right)^{\frac1{p_*}}\\
&\le\gm\pro_{j=1}^N\left(\int_{Q_n}\big|\big[
\ubknp\z_n^{\frac1{p_j}}\big]_{x_j}\big|^{p_j}dx\right)^{\frac1{N p_j}}\\
&\le\gm\pro_{j=1}^N
\left(\bsu_{\ell=1}^N
\frac{2^{p_\ell n}}{\rho_\ell^{p_\ell}}\int_{Q_n}\ubknp^{p_\ell}dx
\right)^{\frac1N\frac1{p_j}}\\
&=\gm\left(\bsu_{\ell=1}^N
\frac{2^{p_\ell n}}{\rho_\ell^{p_\ell}}\int_{Q_n}\ubknp^{p_\ell}dx
\right)^{\frac1{\bp}}.
\end{aligned}
\end{equation*}
From this homogenizing with respect to the measure of $Q_n$ and 
with respect to the integrand, 
\begin{equation}\label{Eq:S:bdd:11}
\begin{aligned}
&\left(\frac1{k^{p_*}}\bint_{Q_{n+1}}\ubknp^{p_*} dx
\right)^{\frac1{p_*}}\\
&\qquad\qquad\le\gm\left(|Q_\rho|^{\frac{\bar p}N}\bsu_{\ell=1}^N
2^{p_\ell n} \frac{k^{p_\ell-\bp}}{\rho_\ell^{p_\ell}} 
\frac{1}{k^{p_\ell}}\bint_{Q_n}\ubknp^{p_\ell}dx
\right)^{\frac1{\bp}}.
\end{aligned}
\end{equation}
For each $\ell\in\{1,\dots,N\}$, estimate
\begin{equation*}
\frac1{k^{p_\ell}}\bint_{Q_n}\ubknp^{p_\ell}dx\le 
\left(\frac1{k^{p_*}}\bint_{Q_n}\ubknp^{p_*}dx
\right)^{\frac{p_\ell}{p_*}}\left(\frac{\big|
[u>\bar{k}_n]\cap Q_n\big|}{|Q_n|}\right)^{1-\frac{p_\ell}{p_*}}.
\end{equation*}
Also
\begin{equation*}
\begin{aligned}
\frac1{k^{p_*}}\bint_{Q_n}\uknp^{p_*}dx&\ge
\frac1{k^{p_*}}\bint_{Q_n\cap[u>\bar{k}_n]}\big(\bar{k}_n-k_n\big)^{p_*}dx\\
&\ge\frac1{2^{p_*(n+2)}}\frac{\big|
[u>\bar{k}_n]\cap Q_n\big|}{|Q_n|}.
\end{aligned}
\end{equation*}
Therefore,
\begin{equation*}
\frac1{k^{p_\ell}}\bint_{Q_n}\ubknp^{p_\ell}dx\le 2^{(p_*-p_\ell)(n+2)} 
\frac1{k^{p_*}}\bint_{Q_n}\uknp^{p_*}dx.
\end{equation*}
Combine these calculations in (\ref{Eq:S:bdd:11}), to get
\begin{equation*}
\begin{aligned}
&\left(\frac1{k^{p_*}}\bint_{Q_{n+1}}(u-k_{n+1})_+^{p_*}dx
\right)^{\frac1{p_*}}\\
&\qquad\qquad \le\gm 2^{\frac{p_*}{\bp}n} 
\left[|Q_\rho|^{\frac{\bar p}N}\bsu_{\ell=1}^N
\frac{k^{p_\ell-\bp}}{\rho_\ell^{p_\ell}}\right]^{\frac1{\bar p}}
\left[\left(\frac1{k^{p_*}}\bint_{Q_{n}}\uknp^{p_*}dx
\right)^{\frac1{p_*}}\right]^{\frac{p_*}{\bp}}.
\end{aligned}
\end{equation*}
Set
\begin{equation}\label{Eq:S:bdd:YN}
Y_n=\left(\frac1{k^{p_*}}\bint_{Q_n}\uknp^{p_*}dx
\right)^{\frac1{p_*}},
\end{equation}
and rewrite the previous inequalities in the form
\begin{equation}\label{Eq:S:bdd:13}
Y_{n+1}\le\gm 2^{\frac{p_*}{\bp}n} 
\left[|Q_\rho|^{\frac{\bar p}N}\bsu_{\ell=1}^N
\frac{k^{p_\ell-\bp}}{\rho_\ell^{p_\ell}}\right]^{\frac1{\bar p}} 
Y_n^{1+\frac{p_*-\bp}{\bp}}.
\end{equation}
Recall that the radii $\rho_j$ are still 
to be chosen. 
\subsection{A Quantitative, Homogeneous Estimate for $\pmx<p_*$}\label{S:bdd:2}
Choose 
\begin{equation}\label{Eq:S:bdd:1}
\rho_j=\rho^{\frac\al{p_j}},
\end{equation}
where $\al$ is an arbitrary positive parameter. Stipulate 
to take $k\ge1$ and estimate
\begin{equation*}
\left[|Q_\rho|^{\frac{\bar p}N}\bsu_{\ell=1}^N
\frac{k^{p_\ell-\bp}}{\rho_\ell^{p_\ell}}\right]^{\frac1{\bar p}}
\le 2 N k^{\frac{\pmx-\bp}{\bp}}. 
\end{equation*}
For such choices (\ref{Eq:S:bdd:13}) yield
\begin{equation}\label{Eq:S:bdd:14}
Y_{n+1}\le\gm 2^{\frac{p_*}{\bp}n}k^{\frac{\pmx-\bp}{\bp}}
Y_n^{1+\frac{p_*-\bp}{\bp}}
\end{equation}
for a new constant $\gm$ depending only on $\{N,p_1,\dots,p_N\}$. 
It follows from these that $\{Y_n\}\to0$ as $n\to\infty$, 
provided
\begin{equation*}
Y_o=\frac1k\left(\bint_{Q_\rho}u_+^{p_*}dx\right)^{\frac1{p_*}}\le 
\gm^{-\frac{\bp}{p_*-\bp}}
2^{-\frac{p_*}{\bp}\left(\frac{\bp}{p_*-\bp}\right)^2} 
k^{-\frac{\pmx-\bp}{p_*-\bp}}. 
\end{equation*}
Thus, choosing 
\begin{equation*}
k=\gm^{\frac{\bp}{p_*-\pmx}}
2^{\frac{p_*}{p_*-\pmx}\frac{\bp}{p_*-\bp}} 
\left[\left(\bint_{Q_\rho}u_+^{p_*}dx\right)^{\frac1{p_*}}
\right]^{\frac{p_*-\bp}{p_*-\pmx}} 
\end{equation*}
yields
\begin{equation}\label{Eq:S:bdd:sup}
\essup_{Q_{\frac12\rho}} u_+\le 1\wedge C 
\left[\left(\bint_{Q_\rho}u_+^{p_*}dx\right)^{\frac1{p_*}}
\right]^{\frac{p_*-\bp}{p_*-\pmx}}, 
\end{equation}
where 
\begin{equation*}
C=\gm^{\frac{\bp}{p_*-\pmx}}
2^{\frac{p_*}{p_*-\pmx}\frac{\bp}{p_*-\bp}}.
\end{equation*}
Write now (\ref{Eq:S:bdd:sup}) over the pair of  cubes 
$Q_{\sig\rho}\subset Q_\rho$, where $\sig\in(\frac12,1)$ is 
an interpolation parameter. Then  
\begin{equation}\label{Eq:S:bdd:15}
\essup_{Q_{\frac12\rho}} u_+\le 1\wedge C^{\prime} 
\left(\bint_{Q_\rho}u_+^{\pmx}dx\right)^{\frac1{\bp}\frac{p_*-\bp}{p_*-\pmx}}.
\end{equation}
\begin{remark}\label{Rmk:S:bdd:1} {\normalfont 
The estimates in (\ref{Eq:S:bdd:sup}) and (\ref{Eq:S:bdd:15}) are 
homogeneous with respect to the cube $Q_\rho$, i.e., they are 
invariant for dilations of the variables $(x_1,\dots,x_N)$ that 
keep invariant the relative intrinsic geometry of 
\eqref{Eq:S:intr:2} and \eqref{Eq:S:bdd:1}. 
In this sense they are an improvement with respect to the 
estimates of Kolod{\={\i}}{\u\i} \cite[Theorem~2]{K}.
If $p_j=\bp$ for all $j=1,\dots,N$ this reproduces the 
classical estimate for isotropic elliptic equations. 
}
\end{remark}
\begin{remark}\label{Rmk:S:bdd:3} {\normalfont 
The constants $C$ and $C^\prime$ in (\ref{Eq:S:bdd:sup}) and 
(\ref{Eq:S:bdd:15}), can be quantitatively determined only in terms 
of $N$ and the $p_j$ for $j=1,\dots,N$. However, they tend to infinity 
as $\pmx\nearrow p_*$. 
}
\end{remark}
\subsection{A Quantitative, Homogeneous Estimate for $\pmx=p_*$}\label{S:bdd:3}
Redefine the levels in (\ref{Eq:S:bdd:7}) as 
\begin{equation}\label{Eq:S:bdd:16}
k_n=\Big(1-\frac1{2^{n+1}}\Big) k,\quad\text{ and }\quad
\bar{k}_n=\frac{k_{n+1}+k_n}2,\qquad\text{ for $n=0,1,\dots$.}
\end{equation}
This implies that $k_o=\frac12k$ and $k_\infty=k$. 
All estimations remain unchanged and yield (\ref{Eq:S:bdd:13}), 
with a slight modification of the constant $\gm$ and with the same 
definition (\ref{Eq:S:bdd:YN}) of the $Y_n$. Continue to choose 
$\rho_j=\rho^{\frac\al{p_j}}$, and stipulate to take $k\ge1$. This yields
the analogues of (\ref{Eq:S:bdd:14}) with $\pmx=p_*$, i.e.,
\begin{equation*}
Y_{n+1}\le\gm 2^{\frac{p_*}{\bp}n}k^{\frac{p_*-\bp}{\bp}}
Y_n^{1+\frac{p_*-\bp}{\bp}}. \tag*{(\ref{Eq:S:bdd:14})${}_{\pmx=p_*}$}
\end{equation*}
Taking into account the definition (\ref{Eq:S:bdd:YN}) 
of the $Y_n$, in this last inequality, the parameter $k$ 
scales out. Thus, setting 
\begin{equation*}
X_n=\left(\bint_{Q_n}\uknp^{p_*}dx\right)^{\frac1{p_*}},
\end{equation*}
the recursive inequalities (\ref{Eq:S:bdd:14})${}_{\pmx=p_*}$ are
\begin{equation}\label{Eq:S:bdd:17}
X_{n+1}\le\gm 2^{\frac{p_*}{\bp}n} X_n^{1+\frac{p_*-\bp}{\bp}}. 
\end{equation}
It follows from these that $\{X_n\}\to0$ as $n\to\infty$, provided
\begin{equation}\label{Eq:S:bdd:18}
X_o= \left(\bint_{Q_\rho}\big(u-{\txty\frac12}k\big)_+^{p_*}
dx\right)^{\frac1{p_*}}\le\gm^{-\frac{\bp}{p_*-\bp}} 
2^{-\frac{p_*}{\bp}\left(\frac{\bp}{p_*-\bp}\right)^2}.
\end{equation}
Since $u\in L^{p_*}_{\loc}(E)$, such a $k$ can be 
quantitatively, although not explicitely, determined, 
in terms of $\|u\|_{L^{p_*}(Q_\rho)}$, and then,
\begin{equation}\label{Eq:S:bdd:19}
\essup_{Q_{\frac12\rho}} u_+\le 1\wedge k. 
\end{equation}
\begin{remark}\label{Rmk:S:bdd:4} {\normalfont 
The estimate in (\ref{Eq:S:bdd:19}) is 
homogeneous with respect to the cube $Q_\rho$, i.e., it is 
invariant for dilations of the variables $(x_1,\dots,x_N)$ that 
keep invariant the relative intrinsic geometry of 
\eqref{Eq:S:intr:2} and \eqref{Eq:S:bdd:1}. 
In this sense, it is an improvement with respect to the 
estimates of Fusco-Sbordone \cite[Theorem~1]{FS}.
}
\end{remark}
\begin{remark}\label{Rmk:S:bdd:5} {\normalfont 
In (\ref{Eq:S:bdd:17}) the number $\kp=\frac{p_*-\bp}{\bp}$ by 
which the power of $X_n$ exceeds one, is precisely determined 
by the estimations, and not arbitrary as it seems to be 
permitted in \cite{FS}.  In view of this, the alternative, 
in the argument of \cite{FS}, by which 
\begin{equation*}
2^{\frac{p_*}{\bp}n} X_n^{1+\frac{p_*-\bp}{\bp}}>1\quad\text{ for 
all $n\ge n_o$  for some $n_o\in\nn$ sufficiently large} 
\end{equation*}
is not needed. Since $n_o$ in \cite{FS} is determined only qualitatively, 
the resulting boundedness estimates seem to be qualitative. 
}
\end{remark}
\begin{remark} {\normalfont 
If one had the additional information that $u\in L^q_{\loc}(E)$, 
for some $q>p_*$, then $k$ in \eqref{Eq:S:bdd:19} could be 
precisely quantified. Indeed, given a non negative 
function $f\in L^q(E)$ and $\varep>0$, consider finding 
$k>0$ such that 
\begin{equation*}
\ine (f-k)^{p}_+dx<\varep \quad\text{ where }\> 0< p<q. 
\end{equation*}
By Chebyshev's inequality $|[f>t]|\le t^{-q}\|f\|_{q,E}^q$, 
for all $t>0$.
Then for $p<q$, 
\begin{equation*}
\begin{aligned}
\ine (f-k)^p_+dx &=p\int_0^\infty s^{p-1} \big|[(f-k)_+>s]\big| ds\\
&=p\int_k^\infty (t-k)^{p-1} \big|[f>t]\big| dt
\le p \int_k^\infty t^{p-1} \big|[f>t]\big| dt\\
&\le  p\,\|f\|_{q,E}^q\int_k^\infty t^{-(q-p)-1} dt
=\frac{p}{q-p}\frac1{k^{q-p}}\|f\|_{q,E}^q. 
\end{aligned}
\end{equation*}
Then choose
\begin{equation*}
k^{q-p}=\frac{p}{\varep}\frac1{q-p}\|f\|_{q,E}^q.
\end{equation*}
}
\end{remark}

\end{document}